\documentclass[11pt,a4paper]{article}

\usepackage[T1]{fontenc}
\usepackage[utf8]{inputenc} 
\usepackage{lmodern}
\usepackage{microtype}

\usepackage[a4paper,margin=1in]{geometry}

\usepackage{amsmath,amssymb,amsthm}

\usepackage{graphicx}
\usepackage{xcolor}
\usepackage{booktabs}

\usepackage[hidelinks]{hyperref}

\usepackage[numbers,sort&compress]{natbib}
\usepackage{authblk}

\theoremstyle{plain}
\newtheorem{theorem}{Theorem}[section]

\newtheorem{corollary}[theorem]{Corollary}
\theoremstyle{definition}

\theoremstyle{remark}

\newcommand{\E}{{\rm E}}
\newcommand{\V}{{\rm Var}}
\newcommand{\Ave}{{\rm Ave}}
\newcommand{\SD}{{\rm SD}}
\newcommand{\PP}{{\rm P}}

\title{Win Rates at First-Passage Times\\ for Biased Simple Random Walks}

\author[1]{\mbox{F.} Thomas Bruss}
\author[1,2]{Davy Paindaveine}

\affil[1]{Universit\'{e} libre de Bruxelles}
\affil[2]{Toulouse School of Economics, Universit\'{e} Toulouse~1 Capitole}

\date{}

\begin{document}

\maketitle

\begin{abstract}
We study the win rate $R_{N_d}/N_d$ of a biased simple random walk~$S_n$ on $\mathbb{Z}$ at the first-passage time $N_d=\inf\{n\ge 0:S_n=d\}$, with~$p=\PP[X_1=+1]\in[1/2,1)$. Using generating-function techniques and integral representations, we derive explicit formulas for the expectation and variance of~$R_{N_d}/N_d$ along with monotonicity properties in the threshold $d$ and the bias $p$. We also provide closed-form expressions and use them to design unbiased coin-flipping estimators of $\pi$ based on first-passage sampling; the resulting schemes illustrate how biasing the coin can dramatically improve \emph{both} approximation accuracy and computational cost.
\end{abstract}


\section{Introduction}
\vspace{2.0001mm}

Let $(S_n)_{n\ge 0}$ be a random walk on $\mathbb{Z}$ with $S_0=0$ and mutually independent increments
$$
X_n:=S_n-S_{n-1}\in\{+1,-1\},
 \quad
\PP[X_n=+1]=p,\ \ \PP[X_n=-1]=q:=1-p,
$$
for~$n\geq 1$. This is one of the most standard settings in the theory of random walks and Markov chains; see, e.g.,
\cite{Durrett2019,LawlerLimic2010,LevinPeres2017}.
Write
\[
R_n:=\#\{1\leq k\leq n:\ X_k=+1\}
\]
for the number of right steps (``wins'') up to time $n$, and, for any positive integer~$d$, define the one-sided hitting time
\[
N_d:=\inf\{n\ge 0:\ S_n=d\}.
\]
Throughout, we restrict to the case~$p\geq 1/2$ for which~$N_d$ is almost surely finite. The object of interest in this work is the \emph{win rate at the stopping time}~$N_d$, namely $R_{N_d}/N_d$, together with its expectation, variance, and related distributional properties. Note that~$R_{N_d}/ N_d$ can also be seen as a performance measure for choosing~$d$ induced by the stopping time~$N_d$.

A first motivation is rooted in goal-based stopping and the reporting bias it induces.
Interpreting each increment as a win/loss outcome, one may view $S_n$ as the net gain after~$n$ rounds and $R_n/n$ as the observed proportion of successes.
In many settings, one does
not fix the number of rounds in advance, but rather stops upon reaching a performance target; here this corresponds to $N_d$, the first time the net gain hits level $d$.
The ratio $R_{N_d}/N_d$ is then exactly the win rate one would report at the moment the target is achieved.
Crucially, this is a \emph{success-conditioned} report: the data are shown only at the time the target is met, which preferentially selects paths with an excess of wins over losses.
Quantifying the resulting bias is a natural way to measure the selection effect created by ``stop upon success'' rules; see, e.g., standard discussions of stopping times and conditioning in random-walk and martingale settings~\cite{AsmussenAlbrecher2010,Williams1991}.

A second motivation is methodological and concerns what classical martingale theory does \emph{not} immediately provide.
Optional stopping theorems control expectations of martingales at stopping times, typically for linear functionals such as $S_{N}$ (and, under additional conditions, for $N$ itself), but they do not directly address nonlinear statistics such as the ratio $R_N/N$.
In particular, even in this elementary random-walk model, the behavior of the empirical proportion $\hat p_N:=R_N/N$ can differ substantially from its fixed-$n$ counterpart once the sample size $N$ is chosen in a path-dependent manner.
Boundary-crossing rules are among the simplest and most classical sequential schemes; see, e.g., \cite{Siegmund1985,TartakovskyNikiforovBasseville2014}.
For the prototypical rule $N_d=\inf\{n\ge 0:\ S_n=d\}$, the statistic $R_{N_d}/N_d$ is precisely the naive empirical proportion evaluated at the first time the cumulative evidence reaches level~$d$.
Thus, quantifying the expectation and variance of~$R_{N_d}/N_d$---and their dependence on $d$ and $p$---provides a concrete, probabilistic answer to a recurring question in sequential procedures: how does optional stopping distort empirical proportions, and how does the distortion scale with the stopping threshold?

The outline of the paper is as follows. In Section~\ref{sec:main}, we derive integral representations for the expectation and variance of the win rate~$R_{N_d}/N_d$ at the first-passage time, and we establish their basic monotonicity and limiting properties as functions of the threshold~$d$ and the bias~$p$. Section~\ref{sec:explicit} then provides closed-form expressions for $\E[R_{N_d}/N_d]$, distinguishing between odd and even values of~$d$. In Section~\ref{sec:pi}, we illustrate how these formulas can be leveraged to construct efficient coin-flipping estimators of~$\pi$ and~$\ln 2$ based on first-passage sampling. Finally, we discuss extensions and limitations in Section~\ref{secConclu}.


\section{Main results}
\label{sec:main}
\vspace{2.0001mm}

Our first main result provides the expectation and the variance of~$R_{N_d}/N_d$. 
\vspace{2.0001mm}

 \begin{theorem}
 \label{TheorMain}
Fix~$p\in[\frac{1}{2},1)$ and a positive integer~$d$. Then,
$$
\E\bigg[\frac{R_{N_d}}{N_d}\bigg]
= \frac12
+
 \frac{d}{2}
\int_0^1
u^{d-1}
h_p(u)
\,
du
$$
and
$$
\V\bigg[ \frac{R_{N_d}}{N_d} \bigg]
=
\frac{d^2}{4}
\bigg\{
\int_0^1 u^{d-1}h_p(u)\ell_p(u)\,du
-
\bigg(
\int_0^1
u^{d-1}
h_p(u)
\,
du
\bigg)^2
\bigg\}
,
$$
where we let
$$
h_p(u)
:=
\frac{p-q u^2}{p+q u^2}
\quad
\textrm{ and }
\quad
\ell_p(u)
:=
\ln\!\Big(\frac{p+qu^2}{u}\Big)
.
$$
 \end{theorem}
\vspace{1.0001mm}

 \begin{proof}
Observe that
$
S_n = R_n - (n - R_n) = 2R_n - n
$
for any $n$. Since $p \ge \tfrac12$, the walk hits $+d$ almost surely, so $N_d $ is finite with probability~$1$. Using~$S_{N_d} = d$, we obtain~$d=2R_{N_d} - N_d$, hence
$$
\frac{R_{N_d}}{N_d}
= \frac12 + \frac{d}{2N_d}
.
$$
Therefore,
\begin{equation}
\label{RelFond}
\E\bigg[ \frac{R_{N_d}}{N_d} \bigg]
=
 \frac12 + \frac{d}{2}
 \,
 \E\bigg[ \frac{1}{N_d} \bigg]
 \quad
 \textrm{ and }
 \quad
 \V\bigg[ \frac{R_{N_d}}{N_d} \bigg]
=
\frac{d^2}{4}
 \,
 \V\bigg[ \frac{1}{N_d} \bigg]
 ,
\end{equation}
and it is sufficient to evaluate the expectation and variance of~$1/N_d$.

Let $T = \min\{ n \ge 0 : S_n = 1\}$ be the first time the walk starting at zero hits level~$+1$. Of course, 
to reach $+d$ from $0$, one must pass through~$
1, 2, \ldots, d
$
in order. By spatial homogeneity and the strong Markov property, the time needed to go from $k$ to $k+1$ has the same distribution as $T$, and these increments are mutually independent. Thus,
\begin{equation}
\label{summ}
N_d \stackrel{\mathcal D}{=} T_1 + T_2 + \cdots + T_d,
\end{equation}
where $T_1, \ldots, T_d$ are mutually independent copies of~$T$. Let the probability generating function of $T$ be~$\phi(t) := \E[t^T]$.
A standard first-step argument yields
$$
\phi(t)
= \frac{1 - \sqrt{1 - 4 p q t^2}}{2qt}
;
$$
see, e.g., (3.6) in Chapter~XI of~\cite{Feller1968} or Section~5.3 in \cite{Zitkovic2019RWLecture5}. Since the $T_i$'s are mutually independent copies of~$T$, the generating function of $N_d$ is then
\begin{equation}
\label{genn}
\phi_d(t) 
:= 
\E[t^{N_d}]
= (\phi(t))^d
= \bigg( \frac{1 - \sqrt{1 - 4 p q t^2}}{2qt} \bigg)^{d}
.
\end{equation}

Using the identity
$$
\frac{1}{n} = \int_0^1 t^{n-1}\,dt, \qquad n \ge 1,
$$
the monotone convergence theorem yields
$$
\E\bigg[\frac{1}{N_d}\bigg]
= \sum_{n=1}^\infty \frac{1}{n}\,\PP[N_d = n] 
= \int_0^1 \sum_{n=1}^\infty \PP[N_d = n]\,t^{n-1}\,dt 
= \int_0^1 \frac{\phi_d(t)}{t}\,dt.
$$
Therefore,
$$
\E\bigg[\frac{1}{N_d}\bigg]
= 
\int_0^1
\frac{1}{t}\bigg( \frac{1 - \sqrt{1 - 4 p q t^2}}{2qt} \bigg)^{d}
\,
dt.
$$
\textcolor{black}{Since~$p\geq q$, we may let
$
t
=
u/(p+q u^2)
$,
which yields}
\begin{equation}
\label{sss}
\E\bigg[\frac{1}{N_d}\bigg]
=
\int_0^1
u^{d-1}
h_p(u)
\,
du
,
\end{equation}
with the function~$h_p$ defined in the statement of the theorem. Plugging this expression in~(\ref{RelFond}) establishes the expectation part of the result. Let us then turn to the variance. The identity
$$
\frac1{n^2}=-\int_0^1 t^{n-1}(\ln t)\,dt,
\qquad
n=1,2,\ldots
$$
(this results from integration by parts)
provides
$$
\E\!\left[\frac1{N_d^2}\right]
=
\sum_{n=1}^\infty \frac1{n^2}\,\PP[N_d=n]
=-\int_0^1 (\ln t)\sum_{n=1}^\infty \PP[N_d=n]\,t^{n-1}\,dt
=-\int_0^1 (\ln t)\,\frac{\phi_d(t)}{t}\,dt
.
$$
Using~(\ref{genn}) and performing the same substitution as above yields
\begin{equation}
\label{zez}
\E\!\left[\frac1{N_d^2}\right]
=
\int_0^1 u^{d-1}h_p(u)\ell_p(u)
\,du
.
\end{equation}
Since~(\ref{RelFond}) provides
$$
\V\bigg[ \frac{R_{N_d}}{N_d} \bigg]
=
\frac{d^2}{4}
 \,
\bigg\{
 \E\bigg[ \frac{1}{N_d^2} \bigg]
-
\bigg(
 \E\bigg[ \frac{1}{N_d} \bigg]
 \bigg)^2
 \bigg\}
,
$$
the result follows from~(\ref{sss}) and~(\ref{zez}).
\end{proof}


Theorem~\ref{TheorMain} allows us to establish the following monotonicity and limiting results.
\vspace{2.0001mm}

 \begin{corollary}
 \label{CorolMonotonicity}
Fix~$p\in[\frac{1}{2},1)$. Then,
$$
\E\bigg[\frac{R_{N_d}}{N_d}\bigg]
>
\E\bigg[\frac{R_{N_{d+1}}}{N_{d+1}}\bigg]
$$
for any positive integer~$d$; moreover,
$$
\lim_{d\to\infty}
\E\bigg[\frac{R_{N_d}}{N_d}\bigg]
=
p
\quad
\textrm{and}
\quad
\lim_{d\to\infty}
\V\bigg[\frac{R_{N_d}}{N_d}\bigg]
=
0
.
$$
 \end{corollary}
\vspace{2.0001mm}

The variance is not monotone in~$d$ in general; for instance, for~$p=9/10$, we have 
$$
\V\bigg[\frac{R_{N_1}}{N_1}\bigg]
<
\V\bigg[\frac{R_{N_2}}{N_2}\bigg]
\quad
\textrm{ and }
\quad
\V\bigg[\frac{R_{N_2}}{N_2}\bigg]
>
\V\bigg[\frac{R_{N_3}}{N_3}\bigg]
.
$$
More importantly, recall that, from a sequential analysis perspective, $R_{N_d}/N_d$ is an estimator of~$p$. Corollary~\ref{CorolMonotonicity} readily entails that this estimator is (positively) biased for any finite~$d$ but that it converges to~$p$ in probability as~$d$ diverges to infinity (since~$0\leq R_{N_d}/N_d\leq 1$, this estimator is therefore asymptotically unbiased). This shows that the distortion introduced by the stopping procedure does not affect consistency.

\begin{proof}
Using Theorem~\ref{TheorMain}, integrating by parts, and noting that~$h_p'(u)<0$ for~$u\in(0,1)$, we obtain
\begin{eqnarray*}
\E\bigg[\frac{R_{N_d}}{N_d}\bigg]
-
\E\bigg[\frac{R_{N_{d+1}}}{N_{d+1}}\bigg]
\!&\!\!=\!\! &\!
 \frac{1}{2}
\int_0^1
\big\{ d\,u^{d-1}-(d+1)u^d \big\}
h_p(u)
\,
du
\\[2mm]
\!&\!\!=\!\! &\!
-
 \frac{1}{2}
\int_0^1
( u^{d}-u^{d+1} )
h_p'(u)
\,
du
>
0
,
\end{eqnarray*}
which establishes the monotonicity result. We turn to the limiting results. Note that the result in Theorem~\ref{TheorMain} rewrites
\begin{equation}
\label{int3}
\E\bigg[\frac{R_{N_d}}{N_d}\bigg]
= 
\frac12
+ \frac{1}{2}
\E[h_p(U_d)]
\end{equation}
and
\begin{equation}
\label{inttttr3}
\V\bigg[ \frac{R_{N_d}}{N_d} \bigg]
=
\frac{1}{4}
\Big\{
d \,
\E[h_p(U_d) \ell_p(U_d)]
-
(\E[h_p(U_d)])^2
\Big\}
,
\end{equation}
\color{black}
where~$U_d$ is such that~$\PP[U_d\leq u]=u^d$ for any~$u\in[0,1]$; equivalently, $U_d$ is
the maximum of~$d$ i.i.d.\ uniform random variables on~$[0,1]$. Since only the marginal
distribution of~$U_d$ matters here, we may realize all of them from a single i.i.d.\
uniform sequence $V_1,V_2,\ldots$ by setting
\begin{equation}
\label{coupling}
U_d:=\max(V_1,\ldots,V_d),
\end{equation}
so that~$U_d\nearrow1$ almost surely. Since~$h_p$ is bounded and continuous, dominated convergence therefore gives
\[
\E[h_p(U_d)]\to h_p(1)=2p-1.
\]
Consequently, (\ref{int3}) yields
\[
\lim_{d\to\infty}
\E\bigg[\frac{R_{N_d}}{N_d}\bigg]
=
\frac12+\frac12h_p(1)
=
p.
\]

Thus, it only remains to prove the variance result, for which a simple proof can also be obtained by exploiting the order-statistic representation of~$U_d$ above. For every positive integer~$d$, we have that~$-\ln U_d=\min(-\ln V_1,\ldots,-\ln V_d)$ is exponential with mean~$1/d$, so that 
\[
U_d\stackrel{\mathcal D}{=}e^{-Z/d}
,
\]
where~$Z$ is exponential with mean one. Since the variance limit depends only on the marginal distribution of~$U_d$,
we may, for the remainder of the proof, use a new coupling based on a single
random variable~$Z\sim\operatorname{Exp}(1)$ and set~$U_d:=e^{-Z/d}$ for
every~$d$. Write~$\alpha:=2p-1$. For every fixed~$z\geq0$, we then have
\begin{equation}
\label{sss1}
h_p(e^{-z/d})\to  h_p(1)=\alpha
\quad \textrm{ as } d\to\infty
\end{equation}
and (from l'H\^{o}pital's rule, after setting~$t=1/d$)
\begin{equation}
\label{sss2}
d\,\ell_p(e^{-z/d})
=
d\ln\big(p+qe^{-2z/d}\big)+z
\to  
\alpha z
\quad \textrm{ as } d\to\infty
.
\end{equation}
Moreover, since
\begin{equation}
\label{NewBounds}
0\leq h_p(u)\leq1
\quad\textrm{ and }\quad
0\leq\ell_p(u)\leq-\ln u
\quad 
\textrm{ for every } u\in(0,1]
\end{equation}
(the result for~$\ell_p(u)$ follows from the identity~$
p+qu^2-u=(1-u)(p-qu)\geq0$ that is valid for~$p\geq 1/2$, together with~$p+qu^2\leq1$), we have
\[
0
\leq
d\,h_p(e^{-Z/d})\ell_p(e^{-Z/d})
\leq
Z
\]
almost surely. From~(\ref{sss1})--(\ref{sss2}), dominated convergence then yields
\[
d\,\E[h_p(U_d)\ell_p(U_d)]
\to  
\alpha^2\E[Z]
=
\alpha^2.
\]
Combining this with~$\E[h_p(U_d)]\to\alpha$ in~(\ref{inttttr3}) proves that
$\V[R_{N_d}/N_d]\to  0$.
\end{proof}


\color{black}
The order-statistic representation of~$U_d$ identified in this proof also gives a direct probabilistic proof of the monotonicity result in Corollary~\ref{CorolMonotonicity}. Indeed, under the coupling in~(\ref{coupling}),
$
U_{d+1}
=
\max(U_d,V_{d+1})
\geq
U_d$
almost surely,
with strict inequality with positive probability. Since~$h_p$ is strictly decreasing on~$[0,1]$, it follows that
$
\E[h_p(U_{d+1})]
<
\E[h_p(U_d)]$, which, in view of~(\ref{int3}),  provides another proof of the claimed monotonicity.
\vspace{2.0001mm}
\color{black}

 \begin{corollary}
 \label{CorolMonotonicityInp}
Fix a positive integer~$d$. Then,
\begin{equation}
\label{s}
p\mapsto \E\bigg[\frac{R_{N_d}}{N_d}\bigg]
\end{equation}
is strictly increasing over $[\frac{1}{2},1]$; moreover,
$$
\lim_{p\to 1}
\E\bigg[\frac{R_{N_d}}{N_d}\bigg]
=
1
\quad
\textrm{ and }
\quad
\lim_{p\to 1}
\V\bigg[\frac{R_{N_d}}{N_d}\bigg]
=
0
.
$$
 \end{corollary}
\vspace{2.0001mm}

Again, the variance of~$R_{N_d}/N_d$ does not show a similar monotonic behavior; in particular, for~$d\geq 3$, the right derivative of~$p\mapsto \V[R_{N_d}/N_d]$ at~$1/2$ is positive, but this variance must \textcolor{black}{decrease for at least some larger values of~$p$,} since it converges to zero as~$p\to 1$.

\begin{proof}
Using Theorem~\ref{TheorMain}, \textcolor{black}{differentiation under the integral sign, justified by dominated convergence,} yields that
$$
\frac{d}{dp}
\E\bigg[\frac{R_{N_d}}{N_d}\bigg]
= 
 \frac{d}{2}
\int_0^1
u^{d-1}
\frac{\partial h_p(u)}{\partial p}
\,
du
=
d
\int_0^1
\frac{u^{d+1}}{(p+qu^2)^2}
\,
du
>
0
$$ 
for any~$p\in[\frac{1}{2},1)$, \textcolor{black}{where, at~$p=\frac{1}{2}$, the derivative is understood as a right derivative. Thus,} the function in~(\ref{s}) is strictly increasing over~$[\frac{1}{2},1)$. Since $R_{N_d}/N_{d}=1$ almost surely for~$p=1$, this extends to~$[\frac{1}{2},1]$. We turn to the proof of the limiting results. \color{black}
For~$p\in[\frac{1}{2},1)$, 
it follows from~(\ref{NewBounds}) that
$|u^{d-1}h_p(u)|
\leq u^{d-1}$
and
$|u^{d-1}h_p(u)\ell_p(u)|
\leq -u^{d-1} \ln u$
for every~$u\in(0,1)$,
\color{black}
where the upper bounds are integrable on $(0,1)$. For any~$u\in(0,1)$, we have $h_p(u)\to1$ and $\ell_p(u)\to \ln(1/u)$ as $p\to1$. By dominated convergence, we thus have
$$
\int_0^1 u^{d-1}h_p(u)\,du
\to
\int_0^1 u^{d-1}\,du
=
\frac1d
$$
and
$$
\int_0^1 u^{d-1}h_p(u)\ell_p(u) \,du 
\to
-\int_0^1 u^{d-1}(\ln u)\,du
=
\frac1{d^2}
.
$$
Plugging this into Theorem~\ref{TheorMain} yields
$$
\lim_{p\to 1}
\E\bigg[\frac{R_{N_d}}{N_d}\bigg]
=
\frac{1}{2}+\frac{d}{2}\times \frac{1}{d}=1
\quad
\textrm{and}
\quad
\lim_{p\to1}\V\bigg[\frac{R_{N_d}}{N_d}\bigg]
=\frac{d^2}{4}\Big(\frac1{d^2}-\frac1{d^2}\Big)=0
,
$$
which concludes the proof.
\end{proof}


\section{An explicit expression for the expectation}
\label{sec:explicit}
\vspace{2.0001mm}

The following result provides an explicit expression for~$\E[R_{N_d}/N_d]$.  
\vspace{2.0001mm}

 \begin{theorem}
 \label{TheorParticularCases}
Fix~$p\in[\frac{1}{2},1)$ and a positive integer~$d$. Then, 
$$
\E\!\left[\frac{R_{N_d}}{N_d}\right]
=
(-1)^{(d-1)/2}
d
\left(\frac{p}{q}\right)^{d/2}
\arctan
\bigg(\sqrt{\frac{q}{p}}\bigg)
+
1
-
d
\sum_{k=0}^{(d-1)/2}
\frac{(-1)^{k}}{d-2k}
\left(\frac{p}{q}\right)^{k}
$$
for $d$ odd, and
$$
\E\!\left[\frac{R_{N_d}}{N_d}\right]
=
(-1)^{d/2}
\frac{d}{2}
\left(\frac{p}{q}\right)^{d/2}
(\ln p) 
+
1
- 
\frac{d}{2}
\left(\frac{p}{q}\right)^{d/2}
\sum_{k=1}^{d/2} 
(-1)^{(d/2)+k}
\binom{d/2}{k}\frac{p^{-k}-1}{k}
$$
for $d$ even.
 \end{theorem}
 \vspace{1mm}

It follows from this result that
$$
\lim_{p\to 1}
\E\bigg[\frac{R_{N_1}}{N_1}\bigg]
=
\lim_{p\to 1}
\sqrt{\frac{p}{q}}
\,
\arctan\bigg(\sqrt{\frac{q}{p}}\bigg)
=
1
,
$$
in accordance with Corollary~\ref{CorolMonotonicityInp}. For~$p = q = \tfrac12$, we obtain for instance the values
\begin{equation}
\label{esttt}
\E\bigg[\frac{R_{N_1}}{N_1}\bigg]
= \frac{\pi}{4}
,
\qquad
\E\bigg[\frac{R_{N_2}}{N_2}\bigg]
=
\ln 2
,
\quad
\textrm{and}
\quad
\E\bigg[\frac{R_{N_3}}{N_3}\bigg]
=
3\Big(1-\frac{\pi}{4}\Big)
.
\end{equation}

 \begin{proof} 
Since
$$
h_p(u) = 1 - \frac{2q u^2}{p+q u^2},
$$
Theorem~\ref{TheorMain} provides
\begin{equation}
\label{ici}
\E\bigg[\frac{R_{N_d}}{N_d}\bigg]
=
\frac12
+ \frac{d}{2}
\bigg(
\frac1d
 - 
 2q\int_0^1 \frac{u^{d+1}}{p+q u^2}\,du
\bigg)
=
1
 - 
 dr\int_0^1 \frac{u^{d+1}}{1+r u^2}\,du
,
\end{equation}
where we let~$r:=q/p$. We treat separately the cases~$d$ odd and~$d$ even.
\vspace{3mm}

(a) Assume that~$d=2m+1$ for some nonnegative integer~$m$. We need to compute 
\begin{equation}
\label{DefFk}
\E\bigg[\frac{R_{N_d}}{N_d}\bigg]
=
1
 - 
 dr
 F_{m+1}
 ,
\quad
\textrm{ with }
\
F_k:=\int_0^1 \frac{u^{2k}}{1+r u^2}\,du
,
\quad
k=0,1,2,\ldots
\end{equation}
For $k\geq 1$,
$$
\frac{u^{2k}}{1+r u^2}
=
\frac{1}{r}
\bigg(u^{2k-2}-\frac{u^{2k-2}}{1+r u^2}\bigg),
$$
so integrating from $0$ to $1$ gives
\begin{equation}
\label{Recursion}
F_k=\frac{1}{r(2k-1)}-\frac{1}{r}F_{k-1}.
\end{equation}
Iterating this recursion from the base case
$
F_0
=
\frac{1}{\sqrt r}\arctan(\sqrt r)
$
yields
$$
F_{m+1}
=
\sum_{j=1}^{m+1}
\frac{(-1)^{j-1}}{2(m+1-j)+1}
r^{-j}
+
\frac{(-1)^{m+1}}{r^{m+1}}F_0
.
$$
Therefore,
$$
\E\bigg[\frac{R_{N_d}}{N_d}\bigg]
=
1
 - 
 d
\bigg(
\sum_{j=1}^{m+1}
\frac{(-1)^{j-1}}{2(m+1-j)+1}
r^{-(j-1)}
+
\frac{(-1)^{m+1}}{r^{m+(1/2)}}
\arctan(\sqrt r)
\bigg)
$$
$$
=
1
 - 
 d
\sum_{k=0}^{m}
\frac{(-1)^k}{2(m-k)+1}
r^{-k}
+
\frac{(-1)^{m}d}{r^{m+(1/2)}}
\arctan(\sqrt r)
,
$$
which shows the result.
\vspace{3mm}

(b) Assume that~$d=2m$ for some positive integer~$m$. Setting
$
u
=
\sqrt{(x-1)/r}
$
in~(\ref{ici}), we need to compute 
$$
\E\bigg[\frac{R_{N_d}}{N_d}\bigg]
=
1
 - 
 dr\int_0^1 \frac{u^{d+1}}{1+r u^2}\,du
=
1
 - 
\frac{d}{2r^m}
 \int_1^{1+r}
 \frac{(x-1)^{m}}{x}\,dx
.
$$
Expanding~$(x-1)^m$, we obtain
\begin{eqnarray*}
\int_{1}^{1+r}\frac{(x-1)^m}{x}\,dx
\!&\!\!=\!\!&\!
\sum_{k=0}^m \binom{m}{k}(-1)^{m-k}\int_{1}^{1+r}x^{k-1}\,dx
\\[2mm]
\!&\!\!=\!\!&\!
(-1)^m\ln(1+r)+\sum_{k=1}^m \binom{m}{k}(-1)^{m-k}\frac{(1+r)^k-1}{k}
.
\end{eqnarray*}
Therefore,
$$
\E\bigg[\frac{R_{N_d}}{N_d}\bigg]
=
1
+ 
(-1)^{m+1}
\frac{d}{2r^m}
\ln(1+r)
 - 
\frac{d}{2r^m}
\sum_{k=1}^m 
(-1)^{m-k}
\binom{m}{k}\frac{(1+r)^k-1}{k}
,
$$
which establishes the result.
\end{proof}

The closed-form expressions in Theorem~\ref{TheorParticularCases} can be numerically delicate when~$d$ is moderate or large, since they involve alternating sums and may therefore suffer from cancellation between terms of comparable magnitude. For stable numerical evaluation, 
\color{black}
it is preferable to evaluate the positive integral in~(\ref{ici}) directly. 
\color{black}


\section{An application: \textcolor{black}{approximating~$\pi$ and~$\ln 2$ with a coin}}
\label{sec:pi}
\vspace{2.0001mm}

\color{black}
The use of random experiments to approximate~$\pi$ is classical; we mention only two approaches. A particularly simple one is to draw a circle inscribed in a square piece of cardboard, expose it briefly to moderate rainfall, and count the raindrops that fall inside the circle. Assuming that the drops are uniformly distributed, the proportion of drops inside the circle relative to the total number falling on the square approximates the area ratio~$\pi/4$.

\color{black}
A more picturesque---and therefore perhaps more widely discussed---example is Buffon's needle problem, although its analysis requires some calculus; see, for example,~\cite{KlainRota1997}. If a needle of length~$L$ is dropped at random onto a plane ruled by parallel lines spaced a distance~$a\geq L$ apart, then the probability that the needle crosses one of the lines is
\[
\frac{2L}{\pi a}.
\]
Repeated needle drops therefore yield a Monte Carlo approximation of~$\pi$. The underlying principle is the same in both experiments.

\begin{figure}[h!]
\centering
\includegraphics[width=0.52\textwidth]{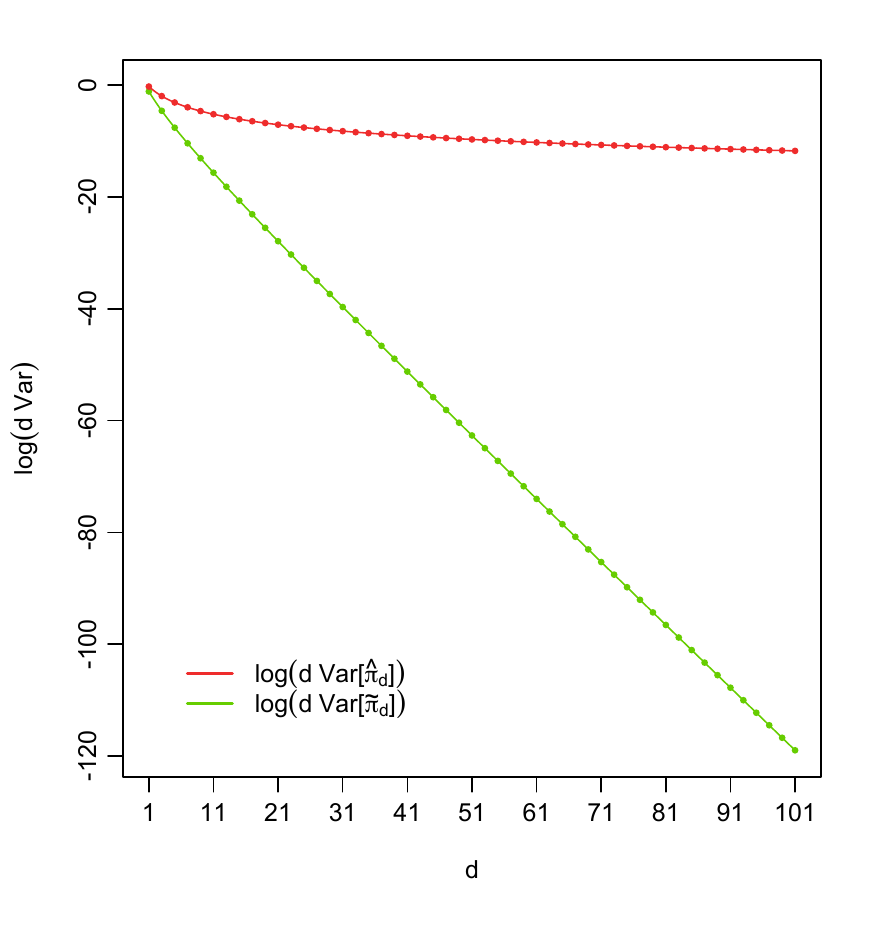}
\vspace{-3mm}
\caption{Plots of $\ln(d\,\V[\hat{\pi}_d])$ and $\ln(d\,\V[\tilde{\pi}_d])$ versus $d\in\{1,3,\ldots,101\}$.
}
\label{Fig1}
\end{figure}

\begin{figure}[h!]
\hspace*{-5mm}
\includegraphics[width=1.05\textwidth]{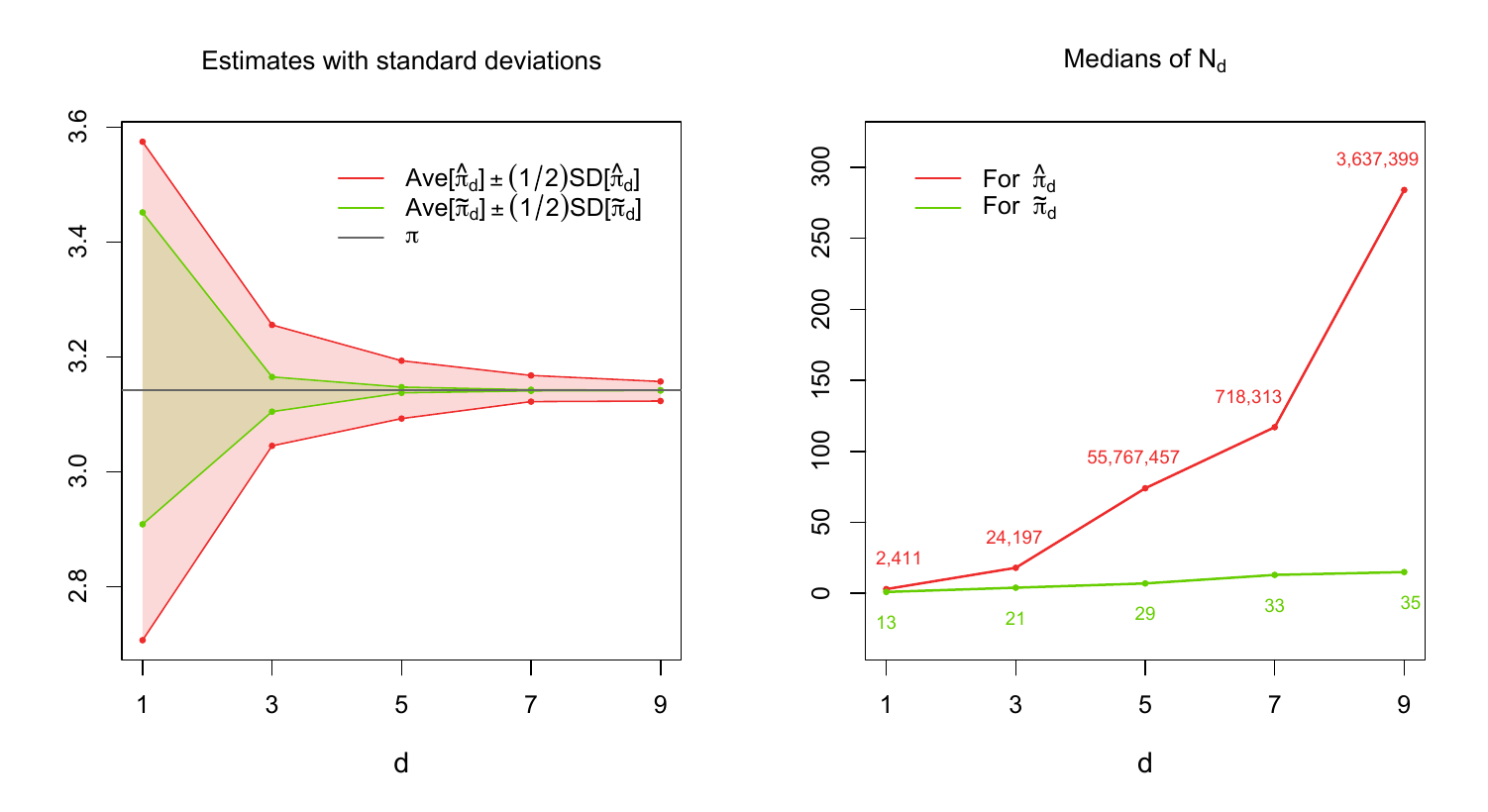}
\vspace{-6mm}
\caption{
Monte Carlo illustrations for the unbiased estimators $\hat{\pi}_d$ (fair coin, $p=\tfrac12$) and~$\tilde{\pi}_d$ (biased coin, $p=\tfrac34$), over odd thresholds $d\in\{1,3,5,7,9\}$.
(Left:) for each $d$, the two shaded bands show the intervals
$\Ave[\hat{\pi}_d]\pm \tfrac12\,\SD[\hat{\pi}_d]$ and
$\Ave[\tilde{\pi}_d]\pm \tfrac12\,\SD[\tilde{\pi}_d]$, based on the averages and standard deviations obtained from  $M=100$ replications; the horizontal line marks $\pi$.
(Right:) sample medians of the corresponding hitting times~$N_d$ (one curve for each estimator), with the numerical labels indicating, for each $d$, the maximum observed value of~$N_d$ among the $M=100$ replications.
}
\label{Fig2}
\end{figure}

In the context of the present paper,
\color{black}
Theorem~\ref{TheorParticularCases} makes it possible to use coin-flipping to approximate irrational constants such as~$\pi$ (for~$d$ odd) and~$\ln 2$ (for~$d$ even). For a fair coin ($p=1/2$), the result in particular implies that
$$
\hat{\pi}_d
:=
(-1)^{(d-1)/2}
\frac{4}{d}
\bigg\{
\frac{R_{N_d}}{N_d}
-
\bigg(
1
-
d
\sum_{k=0}^{(d-1)/2}
\frac{(-1)^k}{d-2k}
\bigg)
\bigg\},
\qquad
d \textrm{ odd},
$$
is an unbiased ``estimator'' of~$4\arctan(\sqrt{q/p})=\pi$ \textcolor{black}{(for~$d=1$, this fair-coin result was obtained independently in~\cite{Propp2026}; see also~\cite{WagenmakersWaldorp2026})}. The approximation error can be controlled by Chebyshev's inequality via
$$
\E[(\hat{\pi}_d-\pi)^2]
=
\V[\hat{\pi}_d]
=
\frac{16}{d^2}
\V\!\left[\frac{R_{N_d}}{N_d}\right]
.
$$
This variance, which can be evaluated by using Theorem~\ref{TheorMain}, converges rather quickly to zero as~$d$ diverges to infinity; see Figure~\ref{Fig1}. However, this approximation strategy is not satisfactory in practice: for~$p=1/2$, the number of flips~$N_d$ is almost surely finite, but its distribution is very heavy-tailed with an infinite expectation, resulting in a computational burden to evaluate~$\hat{\pi}_d$ that quickly deteriorates with~$d$; see the right panel of Figure~\ref{Fig2}. This panel shows that the medians of~$N_d$ obtained from 100 independent replications increase faster than linearly with~$d$ (we consider the median since~$\E[N_d]$ is infinite) and that the maximal value of~$N_d$ over these replications may be huge.

This can be improved significantly by using a biased coin. Indeed,  the expected computational cost~$\E[N_d]=d/(2p-1)$ is then finite and increases only linearly in~$d$ (obviously, this cost  decreases with~$p$). Since~$6\arctan(\sqrt{1/3})=\pi$, we take~$p=3/4$ and consider the unbiased estimator
$$
\tilde{\pi}_d
:=
(-1)^{(d-1)/2}
\frac{6}{3^{d/2}d}
\bigg\{
\frac{R_{N_d}}{N_d}
-
\bigg(
1
-
d
\sum_{k=0}^{(d-1)/2}
\frac{(-3)^k}{d-2k}
\bigg)
\bigg\}
,
\qquad
d \textrm{ odd},
$$
of~$\pi$ resulting from Theorem~\ref{TheorParticularCases}. When investigating the dependence on~$d$ of the corresponding approximation error, it is natural to consider the quantity
\begin{equation}
\label{faircomp}
d\,\V[\tilde{\pi}_d]
=
\frac{4}{3^{d-2} d}
\V\!\left[\frac{R_{N_d}}{N_d}\right]
;
\end{equation}
indeed, since the computational cost increases linearly with~$d$, it is equally costly on average to evaluate~$\tilde{\pi}_d$ as to obtain~$d$ mutually independent estimators~$\tilde{\pi}_1^{(1)},\ldots,\tilde{\pi}_1^{(d)}$, which would provide the averaged estimator~$\frac{1}{d}\sum_{k=1}^d \tilde{\pi}_1^{(k)}$ that has variance~$
\frac{1}{d} \V[\tilde{\pi}_1]
$.
Therefore, the quantity in~\eqref{faircomp} gives a way to study how the approximation error depends on $d$ while keeping the computational budget fixed. Figure~\ref{Fig1} indicates that this error proxy decays exponentially fast to~$0$. As a result, for sufficiently large~$d$, the accuracy gain from increasing~$d$ dominates the linear growth of the computational cost by any prescribed factor. This behavior is corroborated in Figure~\ref{Fig2}, where the biased-coin procedure is seen to outperform the fair-coin one  both in terms of approximation accuracy and computational cost.

In practice, one needs to choose a value of~$d$ that will ensure a level of accuracy with \textcolor{black}{a prescribed probabilistic guarantee}. 
Applying Chebyshev's inequality, we obtain
$$
\PP[
|\tilde{\pi}_d-\pi|
>
\varepsilon]
\leq
\frac{1}{\varepsilon^2}
\V[\tilde{\pi}_d]
=
\frac{4}{3^{d-2} d^2\varepsilon^2}
\V\!\left[\frac{R_{N_d}}{N_d}\right]
,
$$
where the variance can be  evaluated by using Theorem~\ref{TheorMain}. For any given~\textcolor{black}{$\delta\in(0,1)$ and~$\varepsilon>0$}, this can be used to identify a value of~$d$ for which~$|\tilde{\pi}_d-\pi|\leq \varepsilon$ with probability at least~$1-\delta$. For instance, with~$d=45$, the probability that~$|\tilde{\pi}_d-\pi|\leq 10^{-9}$ is at least~$1-10^{-6}$. While~$d=45$ may seem large, the bias of the considered random walk keeps the computational burden under control. In~$10^4$ independent replications with~$d=45$, the maximal value of~$|\tilde{\pi}_d-\pi|$ was below~$10^{-12}$, and the minimal, \textcolor{black}{mean}, and maximal values of~$N_d$ were~47, 89.9, and~181, respectively (for~$p=3/4$,~$N_d$ has expectation~$2d$ and standard deviation~$\sqrt{6d}$).

It is natural to wonder whether Theorem~\ref{TheorParticularCases} may lead to an even better procedure to approximate~$\pi$.
Obviously, the larger~$p$, the smaller the computational cost. 
Moreover, for odd~$d$, the leading term in Theorem~\ref{TheorParticularCases} involves the factor $(p/q)^{d/2}$, which suggests that increasing~$p$ should also improve accuracy at a fast rate. Two remarks are in order.

\begin{itemize}
\item
If one insists that the ratio $r:=q/p$ be \emph{rational} and that $\arctan(\sqrt r)$ be a rational multiple of~$\pi$, then there is no gain beyond $r=1/3$.
Indeed, if $\theta=\arctan(\sqrt r)=m\pi/n$ with  positive integers~$m,n$, then
\[
\cos(2\theta)=\frac{1-\tan^2\theta}{1+\tan^2\theta}=\frac{1-r}{1+r}
\]
is rational. Since $2\theta\in\pi\mathbb{Q}$, \textcolor{black}{the theorem} in~\cite{Olmsted1945RationalTrig}
(see also Corollary~3.12 in~\cite{Niven1956IrrationalNumbers})
implies that $\cos(2\theta)\in\{0,\pm\tfrac12,\pm 1\}$, hence $r\in\{\tfrac13,1,3\}$.
In particular, among rational $r>0$ such that $\arctan(\sqrt r)\in\pi\mathbb{Q}$, the
smallest possible value is $r=1/3$.
\vspace{2mm}
\item
The situation is different if we allow $r$ to be an algebraic real number, though. For any integer $k\ge 4$, set
\[
\theta_k:=\frac{\pi}{k},
\qquad
r_k:=\tan^2(\theta_k),
\qquad
p_k:=\frac{1}{1+r_k},
\quad
\textrm{and}
\quad
q_k:=\frac{r_k}{1+r_k}.
\]
Then, $r_k$ is algebraic (e.g., 
$r_8=(\sqrt2-1)^2$ and $r_{12}=(2-\sqrt3)^2$)
and satisfies $\arctan(\sqrt{r_k})=\theta_k=\pi/k$, so that $\pi=k\,\arctan(\sqrt{r_k})$. Obviously, $r_k\to 0$ as $k\to\infty$.
Thus, one can take $p_k\to 1$ while keeping an exact identity $\pi=k\,\arctan(\sqrt{r_k})$,
thereby reducing the expected cost $\E[N_d]=d/(2p_k-1)$ essentially down to~$d$. For given~$k$, the unbiased estimator of~$\pi$ resulting from Theorem~\ref{TheorParticularCases} is
\[
\tilde{\pi}_{d,k}
:=
(-1)^{(d-1)/2}\,
\frac{kr_k^{d/2}}{d}
\bigg\{
\frac{R_{N_d}}{N_d}
-
\bigg(
1
-
d
\sum_{j=0}^{(d-1)/2}
\frac{(-1)^j}{d-2j}\,r_k^{-j}
\bigg)
\bigg\},
\qquad d\ \text{odd}.
\]
For~$k=6$, this recovers $\tilde{\pi}_d$,
while any $k>6$ would improve on computational cost and approximation accuracy, at the expense of a slightly more complex estimator.
\end{itemize}

The situation is similar when estimating~$\ln 2$ based on the result for even~$d$ in Theorem~\ref{TheorParticularCases}.

\color{black}
To put these stochastic constructions in perspective, classical deterministic methods provide rational approximations of~$\pi$ and~$\ln 2$ with explicit error guarantees. For instance, Machin's identity 
$
\tfrac{\pi}{4}
=
4\arctan(\tfrac{1}{5})
-
\arctan(\tfrac{1}{239})
$
(see, e.g., (6) in~\cite{BorweinBorweinGalway2004})
and the power-series expansions of the arctangent and logarithm lead to the rational approximations
$$
\pi_m
:=
16\sum_{k=0}^{m}
\frac{(-1)^k}{(2k+1)5^{2k+1}}
-
4\sum_{k=0}^{m}
\frac{(-1)^k}{(2k+1)239^{2k+1}}
\quad\textrm{ and }\quad
\lambda_m
:=
\sum_{k=1}^{m}\frac{1}{k2^k}
$$
of~$\pi$ and~$\ln 2$, respectively. The corresponding truncation errors satisfy
$$
|\pi-\pi_m|
\leq
\frac{16}{(2m+3)5^{2m+3}}
+
\frac{4}{(2m+3)239^{2m+3}}
\quad\textrm{ and }\quad
0<(\ln 2)-\lambda_m
\leq
\frac{1}{(m+1)2^m}.
$$
Thus, these deterministic errors are nonrandom and decrease geometrically with~$m$,
whereas the errors of our unbiased estimators are random and are assessed through their
variance relative to the expected computational cost. Concretely, to achieve an absolute error of at most~$10^{-9}$, $m=6$ suffices for~$\pi_m$ (that is, 14~rational summands), since
$
|\pi-\pi_6|\leq 3.5\times10^{-11}
$
(for~$\ln 2$, the displayed error bound requires~$m=26$). By comparison, with~$p=3/4$, we saw above that~$\tilde{\pi}_{45}$ requires~$45/(2p-1)=90$ coin flips on average, and Chebyshev's inequality ensures that
$
|\tilde{\pi}_{45}-\pi|\leq10^{-9}
$
with probability at least~$1-10^{-6}$, rather than with certainty (of course, a summand in Machin's formula and a coin flip are not directly comparable elementary operations).
Deterministic methods are therefore preferable when the sole objective is high-precision numerical computation, even though the standard deviation of~$\tilde{\pi}_d$ also decreases at least geometrically in~$d$, since~(\ref{faircomp}) and~$0\leq R_{N_d}/N_d\leq 1$ yield
\[
\SD[\tilde{\pi}_d]
=
\frac{6}{3^{d/2}d}
\sqrt{\V\!\left[\frac{R_{N_d}}{N_d}\right]}
\leq
\frac{3}{3^{d/2}d}.
\] 
The purpose of the present constructions is instead to provide first-passage representations of these constants and to study how the coin bias affects their accuracy--cost trade-off.
\color{black}


\section{Extensions and limitations}
\label{secConclu}

Although we did not address optimality questions in this paper, our results connect naturally to them. For example, fix $\textcolor{black}{\tfrac12}<p<1$ and define the performance of choosing a threshold $d$ as~$\E[R_{N_d}/N_d]$. If $d$ must be selected from a finite set $\{d_1, \ldots, d_s\}$ and picking $d_k$ carries a price~$\pi_k$, then maximizing the price-weighted performance  means identifying the optimal choice of $d$. Similar questions can be extended from random walks to other base processes (e.g., Brownian motion or Poisson jump processes), since random walks are building blocks for many such models. 

One can also consider other directions such as multiple stopping-time problems. For instance, let~$X_1, \ldots, X_n$ be \mbox{i.i.d.} random variables on~$[0,1]$ that are observed sequentially, and assume that one wants to maximize the probability of identifying online, and without recall, the last $k$ records in $X_1, \ldots, X_n$. The optimal strategy was found by~\cite{BrussPaindaveine2000Selecting} in terms of a  generalized odds algorithm. If we define different rewards for $n\ge k$ (for $k$ fixed) as a function of $n$ and achieving the objective, then the optimal reward-weighted performance identifies the optimal strategy at the same time.

Returning to the stopped random walk, the setting becomes more challenging if we allow for 
path-dependent online intervention by changing, if desired, from a  currently chosen $d$ to another $d$. This version becomes interesting if we fix an upper bound $T$ (positive integer) for the available time after which we do not receive a reward. We look, as before, at the list of price tags and choose a starting $d$ as if there was 
 no time constraint, i.e., $ T= \infty$.  Now, if the choice of $d$ can be made at any time, then maximizing the expected value becomes a multistage decision problem
(equivalent to a  problem of optimal control). A strategy is now a sequence of instructions to keep or change the current $d$, still with the same list of possible $d$ and the same price tags. Let~$D$ be the {\it last} (path-dependent) choice of~$d$ before the end of horizon~$T$ according to the sequence of instructions. The goal is now to maximize $\E[R_{N_D}/N_D]$. This problem is computationally much more involved, of course, but clearly it raises no new difficulties because, at each decision time, all relevant information is known: the current level of the random walk, the remaining time, and the additional height to achieve within that remaining time. Consequently, one can evaluate the expected 
final reward for any strategy, i.e., for any sequence of instructions to deal with the current observation.
Since the list of possible choices of $d$ is finite (at most $s$), there exist at most $T\times s \times
d_{\max}$ possible relevant strategies, and, for each of them, we can compute the performance.
 The optimal performance would thus again, like in our first example, identify the optimal strategy.

It would be interesting to go further: ideally, a suitable notion of performance would not only single out an optimal strategy (as in the examples above) but would also help to construct it. Alas, for this objective we encounter limitations, because in some problems we do not see a performance measure that offers a way around the complexity of the original problem. 

Consider, for instance, the celebrated unsolved Robbins' problem of minimizing the expected rank (see, e.g., \cite{BrussFerguson1993}; for a review, see~\cite{Bruss2005}). The problem is the following.

Let the random variables~$X_1, \ldots, X_n$ be \mbox{i.i.d.} uniform on $[0,1]$. We say that~$X_i$ has rank~$k$ if it is the~$k$th smallest value among~$X_1,\ldots,X_n$ (for instance, $X_i$ has rank~$1$ if and only if~$X_i=\min(X_1,\ldots,X_n)$). A sequential observer of   $X_1, \ldots, X_n$, who must stop on exactly one of the $X_k$'s without recall, would like to minimize the expected rank of their choice. Note that if one stops at~$X_k$, its expected rank is $r(X_k)+(n-k)X_k$, where~$r(X_k)$ is the relative rank of~$X_k$ observed so far. This expected rank must be compared with the optimal value by continuing. The latter is, however, fully history-dependent, i.e., depends on all $X_1, \ldots, X_k$; consequently, no simpler sufficient statistic exists. The problem is thus intractable for larger $n$ and the solution is only known for~$n=1,2,3,4$; see~\cite{DendievelSwan2016RobbinsN4}. Note that those of the future values  $X_{k+1}, \ldots, X_n$ which fall into $[0,X_k]$ are again \mbox{i.i.d.} uniform on~$[0,X_k]$, and the problem is essentially the same as in the beginning. We say `essentially' because we are now at time $k$.  This raises the question: which performance measure could potentially reduce the complexity of Robbins' problem? If we could identify one that converges in time sufficiently fast before the end of the horizon as~$n$ grows, this could be crucial.  

We have  to leave our excursion here because we have made no progress worth mentioning on this problem, but not without mentioning that many colleagues have worked on Robbins’ problem, using interesting approaches such as those in~\cite{BrussSwan2009}, \cite{Gnedin2007}, \cite{GnedinIksanov2011}, and~\cite{MeierSogner2017}.  We should also mention that the rich theory linking optimal stopping and free boundary problems~\cite{PeskirShiryaev2006}   motivated the authors of~\cite{BrussSwan2009} to embed Robbins’ problem into continuous time and to present it in terms of a  differential equation. \textcolor{black}{This approach was successful insofar as the continuous-time embedding preserved the limiting value~$v$ and allowed it to be characterized by a single differential equation. It therefore provided a compact analytical characterization of~$v$. The result was nevertheless only a partial success, since the differential equation still involved an auxiliary unknown function, preventing the precise computation of~$v$.} Finally, we should recall that another celebrated stopping problem is linked with the name of Robbins, namely the so-called $S_n/n$-problem; see, e.g., \cite{ChristensenFischer2022}. This problem is related to our (time-dependent) random walk problem introduced above, but more complex and also not completely solved.


\section*{Acknowledgments}

Davy Paindaveine is a member of the Royal Academy of Science, Letters and Fine Arts of Belgium. Thomas Bruss acknowledges the support through his invited professorship at ULB. 


\section*{Funding}

Davy Paindaveine was supported by the “Projet de Recherche” T.0230.24 from
the FNRS (Fonds National pour la Recherche Scientifique), Communauté Française de Belgique.


\bibliography{Arxiv}

\end{document}